\documentclass{article}
\usepackage{amsmath,amssymb,graphicx} %load extra symbols and environments
\usepackage{amssymb}
\usepackage{amsthm}
\usepackage[margin=0.9in]{geometry} %set margins

\usepackage{hyperref}
\usepackage{{tcolorbox}} 
\usepackage{tikz}
\allowdisplaybreaks
\usepackage{comment} 
\usepackage{cancel}
\usepackage{enumitem}  
\usepackage{appendix}
\usepackage{lineno}
\usepackage{xcolor}

\begin{document}
%----------------------------------------
%\nocite{*} % this command forces all references in template.bib to be printed in the bibliography

\title{Global random walk for one-dimensional one-phase Stefan-type moving-boundary problems: Simulation results}

\author{Nicolae Suciu$^{a}$, Surendra Nepal$^{b}$, Yosief Wondmagegne$^{b}$, Magnus Ögren$^{c,d}$,  Adrian Muntean$^{b}$\footnote{Corresponding author, email: adrian.muntean@kau.se}\\
$^{a}$ Tiberiu Popoviciu Institute of Numerical Analysis, Romanian Academy\\
 F\^{a}nt\^{a}nele str.
no. 57, Cluj-Napoca, Romania\\
$^{b}$ Department of Mathematics and Computer Science, Karlstad University,\\
Universitetsgatan 2, Karlstad, 65188, Sweden\\
$^{c}$ School of Science and Technology, Örebro University\\
SE-70182, Örebro, Sweden \\
$^{d}$ HMU Research Center, Institute of Emerging Technologies,\\
GR-71004, Heraklion, Greece \\
\date{\today}}
\maketitle

\noindent

%----------------------------------------
\begin{abstract}
This work presents global random walk approximations of solutions to one-dimensional Stefan-type moving-boundary problems. We are particularly interested in the case when the moving boundary is driven by an explicit representation of its speed. This situation is usually referred to in the literature as moving-boundary problem with kinetic condition.  As a direct application, we propose a numerical scheme to forecast the penetration of small diffusants into a rubber-based material. To check the quality of our results, we compare the numerical results obtained by global random walks either using the analytical solution to selected benchmark cases or relying on finite element approximations with {\em a priori} known convergence rates. It turns out that the global random walk concept can be used to produce good quality approximations of the weak solutions to the target class of problems. 

%\red{previous version: /Dropbox/AdrianMuntean/Paper/version\_3-24.07.2024/GRWM\_Stefan\_Problem (Copy).zip}
	
	
\end{abstract}
%----------------------------------------

%----------------------------------------
\section{Introduction}
\label{sec:intro} 
Models for the penetration of diffusants in rubbers (or in other hyperelastic materials) can under certain conditions be formulated as one-dimensional moving-boundary problems. They are somewhat similar to the classical one-phase Stefan problem describing the melting of ice blocks (see e.g., \cite{Gupta2003, Stefan1891, tarzia1997one}). Briefly speaking, for any pair of space and time variables $(x,t)$, the evolution of the observable $m(x,t)$, denoting the mass concentration of the chemical species of interest,  is governed in our context by the diffusion equation
\[
\frac{\partial m}{\partial t} -\alpha\frac{\partial^2 m}{\partial x^2} = 0 \mbox{ in } Q_s(T),
\]
where $\alpha>0$ is a constant diffusion coefficient. 
This evolution is posed in the % a non-cylindrical 
space-time domain 

$$ Q_s(T):= \{ (x, t) \; | \; x \in (0, s(t))  \; \text{and}\; t \in (0, T)\},$$ where $s(t)$ denotes the maximum depth where the concentration has reached into the material and $T>0$ is a  given final time of the overall process. We refer to this location $s(t)$ as penetration depth and observe that it is not known {\em a priori}.
To fix ideas, we take a  Dirichlet boundary condition at the left (fixed) boundary and the Stefan condition at the right (moving) boundary %(Sect.~\ref{sec:Stefan})
in the classical Stefan problem~\cite{Stefan1891}. The latter essentially relates the speed $s'(t)$ of the moving boundary  to the flux of the transported quantity $m(x,t)$, i.e. 
\[
s^{\prime}(t)\propto -\frac{\partial m}{\partial x}(s(t),t).
\]
In this case, we are led to the classical Stefan condition -- the proportionality constant can be identified precisely and leads to the interface balance law

\[
L\, s^{\prime}(t)=-\alpha \frac{\partial m}{\partial x}(s(t),t),
\]
where $L\geq 0$ is here the corresponding concept of latent heat needed for a purely diffusive setting. 
At this point, it is relevant to note that the classical setting of the Stefan problem takes as boundary condition at the moving boundary $x=s(t)$ the expression $m(s(t),t)=0$. Now, prescribing as well the initial position $s(0)$ of the moving boundary and the initial concentration profile $m(\cdot,t=0)$, the complete formulation of the model is able to describe a wide class of phase transitions scenarios of first-order type; see  \cite{Gupta2003}. 

Within the present framework, we shall consider different information linked to the evolution of the moving boundary, namely we take   %(Sect.~\ref{sec:Stefan_Kinetic}),
\[
s^{\prime}(t)\propto f\left(s(t),m(s(t),t)\right).
\] 
Such a relation is usually referred in the literature as a {\em kinetic condition}; see e.g. \cite{Visintin}. The original interface condition $m(s(t),t)=0$ posed by J. Stefan in \cite{Stefan1891}  to go along with the balance in interface fluxes describes an equilibrium picture. The kinetic representation of the speed describes a non-equilibrium scenario. For instance, if we think of small diffusant particles penetrating rubbers (see Sect.~\ref{sec:rubber}), the mentioned kinetic condition  takes the particular form,
$$
s^{\prime}(t)\propto \left(m(s(t),t) - g(s(t))\right),
$$
where the term $g(s(t))$ mimics the swelling behavior of the rubber. 
To fit well to physically motivated scenarios, the homogeneous Dirichlet boundary condition at the fixed boundary is replaced by a flux boundary condition (of Robin type).
It is worth also noting that the classical one-phase Stefan problem (so the equilibrium picture) admits probabilistic interpretations (see, for instance, \cite{Tsunoda} and related work)-- this gives trust in the quality of the model. On the other hand, in spite that most kinetic formulations are well-posed in suitable function spaces, the one-phase Stefan-type problem with kinetic condition does not seem to have available any obvious probabilistic interpretation. We do not attempt here to unveil such probabilistic interpretation valid for a non-equilibrium phase transition. Instead, we wish to see to which extent a probabilistic-type numerical method can be designed to suit such class of problems.
%The main question we are addressing in this paper is:
This leads to the main question that we are addressing in this paper, namely
 \begin{tcolorbox}
(Q) To which extent a global random walk approach is able to approximate the weak solution of one-phase Stefan-type problems with kinetic condition.  
 \end{tcolorbox}
Looking at (Q), one may wonder from where such sudden interest raises, as in fact one currently knows very well how to approximate numerically  (in a both rigorous and efficient manner) the weak solution to the underlying moving-boundary problem; see e.g. \cite{nepal2023analysis} for details  on how finite element approximations work in this case (and other established numerical methods work as well). Our motivation really stems from the fact that in the engineering application that we have in mind [see the problem description in section \ref{sec:rubber} and compare as well \cite{nepal2021moving} for more background on its potential technologic relevance\footnote{Moving-boundary problems with kinetic boundary conditions have turned to be useful  tools for investigating a number of scenarios (unrelated directly to the ones studied here) for instance in chemical reactor engineering \cite{Villa}, corrosion of cementitious materials \cite{Aiki}, silicon oxidation \cite{Evans}, and much more.}] -- only a finite number of particles succeed in fact to penetrate the hyperelastic material and affect the length of its lifetime. This makes us wonder whether a random walk-type approach could describe perhaps better the physical picture. Such an approach would potentially be able to give insights on how many ingressing particles are needed to damage the host material. In this preliminary study, we wish to address the issue numerically. We proceed in the same spirit as in our previous work \cite{Nepaletal2023}, where the used  algorithms did employ classical random walks, but now we wish to benefit of the global random walk (GRW) structure which seems to lead to much faster calculations.  Recently, the use of GRW was very successful in a number of situations involving reactive transport posed in fixed porous media domains, see, for instance, \cite{Suciuetal2021} and references cited therein. We aim to extend the GRW approach to the case of moving one-dimensional domains. 

This paper is organized in the following fashion:
%\tableofcontents
Section~\ref{sec:grw} provides a brief description of the GRW method as applied for a diffusion problem.  
In the first part of Section~\ref{sec:GRWStefan}, we introduce the one-dimensional classical Stefan problem and its GRW approximation. Additionally, we compare the GRW simulation results with finite element simulation results (with known convergence properties), and finally, we estimate the order of convergence of the GRW scheme. We are mainly interested in approximating moving interfaces whose speed is in some way controlled via a kinetic condition. To this end, we present in the second part of Section~\ref{sec:GRWStefan}, the GRW simulation results for the Stefan-type problem with kinetic condition, where we numerically determine the order of convergence.
%We introduce the one-dimensional classical Stefan problem and its GRW approximation in Section \ref{sec:Stefan}. Additionally, we compare the GRW simulation results with finite element simulation results (with known convergence properties), and finally, we estimate the order of convergence of the GRW scheme. We are mainly interested in approximating moving interfaces whose speed is in some way controlled via a kinetic condition. To this end, in Section  \ref{sec:Stefan_Kinetic}, we present the GRW simulation results for the Stefan-type problem with kinetic condition, where we numerically determine the order of convergence. 
The main aim of this paper is to apply the GRW method to approximate correctly the solution to the moving boundary problem describing the penetration of diffusants into rubber, which is discussed in Section \ref{sec:rubber}. Our conclusion on the obtained results and potential directions for future research are presented in Section \ref{discussion}. 

\section{GRW approximation of diffusion} % Global random walk (GRW) 
\label{sec:grw}
\bigskip
%\red{This section will only present GRW approximation of diffusion; Nicolae's work must be cited here... }

The  theory behind GRW is described in \cite{Suciu2019} and references cited therein. Here we only briefly recall the main basic ingredients. Numerical approximations reported here will be provided by a GRW algorithm, free of numerical diffusion, which approximates the transported quantities with the same precision as that of the finite difference scheme, while keeping the corpuscular description of the corresponding continuous fields \cite[Sect. 3.3]{Suciu2019}. The GRW scheme can be derived as a randomization of the forward-time central-space finite difference  discretization of the diffusion equation,
\[
\frac{m_{i, k+1}-m_{i, k}}{\Delta t}=\alpha\frac{1}{\Delta x}\left(\frac{m_{i+1, k}-m_{i, k}}{\Delta x}-\frac{m_{i,k}-m_{i-1,k}}{\Delta x}\right),
\]
where $m_{i,k}=m({i\Delta x,k\Delta t})$ is the approximation of the continuous concentration field $m(x, t)$ at the lattice sites $i=1,\ldots ,L/\Delta x$ and time points $k=1,\ldots ,T/\Delta t$, $[0,L]\times[0,T]$ is the space-time domain of the problem, and $\Delta x$ and $\Delta t$ are the space and time steps, respectively. The finite difference approximation is represented through the distribution of $\mathcal{N}$ random walkers on the regular lattice,
$m_{i,k}\propto n_{i,k}$. With this, the finite differences scheme becomes
\[
n_{i,k+1}=(1-r)n_{i,k}+\frac{r}{2}(n_{i-1,k}+n_{i+1,k}),
\]
where $r := 2\alpha\,\Delta t/(\Delta x)^2$. 
%$\big( 2\alpha\,\Delta t/{\color{red} (\Delta x)^2} \big)$. 
Accordingly, each group of $n_{i,k}$ particles is distributed over the lattice sites as
\[
n_{j,k}=\delta n_{j,j,k}+\delta n_{j-1,j,k}+\delta n_{j+1,j,k},
\]
where $\delta n_{j,j,k}$ represents the number of particles staying at the initial lattice site $j$, $\delta n_{j-1,j,k}$ is the number of particles moving to the left, and finally, $\delta n_{j+1,j,k}$ is the number of particles moving to the right. If the particles move according to the random walk rule, then $r/2$ corresponds to the jump probability, while the quantities $\delta n$ are binomially distributed random variables. For consistency with the finite differences scheme, their ensemble average, denoted by an overbar in the following expressions, must verify the relations
\[
\overline{\delta n_{j,j,k}}=(1-r)\;\overline{n_{j,k}},\;\; \overline{\delta n_{j\mp 1,j,k}}=\frac{r}{2}\;\overline{n_{j,k}}.
\]
In the GRW approach, instead of generating trajectories for each particle, the solution is obtained by computing the binomial random variables $\delta n$. The latter are accurately approximated with the ``reduced fluctuations algorithm'' \cite[Sect. 3.3.2.2]{Suciu2019} in which the numbers $(1-r)n_{j,k}$ and $rn_{j,k}/2$ are rounded to integers with the {\it floor} function, the fractional parts are summed up, and a new particle is assigned to the lattice site where the sum reaches the unity \cite[Sect.~3.1.1]{SuciuandRadu2022}.

%\begin{remark}\label{rem:GRWconv}
One notices that the normalization to the unity of the random walk probabilities imposes the constraint $r\leq 1$, equivalent to 
%$\big(  \alpha\,\Delta t/{\color{red} (\Delta x)^2} \big) $
$\alpha\,\Delta t/(\Delta x)^2 \leq 1/2$. This is the von Neumann condition that ensures the stability and the convergence of the finite differences scheme considered here \cite[Sect.~6.3]{Strikwerda2004}. On the other side, the relations verified by the quantities $\overline{\delta n}$ imply the equivalence between the ensemble-averaged GRW and the finite differences approximations.  Note that the GRW algorithm is self-averaging in the sense that the single realization solution $n_{j,k}$ approaches its ensemble average $\overline{n_{j,k}}$ with the increase of the number of particles; see  \cite[Sect.~3.3.2.3]{Suciu2019}. For large enough $\mathcal{N}$, the GRW and finite differences approximations are therefore identical in the limit of the machine precision, which implies the expected convergence of the GRW scheme. See also a numerical demonstration of the self-averaging of the GRW scheme for a non-linear diffusion equation (Richards' equation) in \cite[Sect. 3.1]{Suciuetal2024}.
%\end{remark}

All the GRW simulations presented in this article use the fixed number of random walkers $\mathcal{N}=10^{24}$. That is, in particular, the concentrations are represented by mole fractions. %\red{new text: $\rightarrow\rightarrow$} 
%The simulations are implemented in MATLAB codes with a generic structure as follows:
The algorithms are implemented in MATLAB with a generic structure as follows: a main code, where the initial data and the boundary conditions are formulated for transported quantities $m(x,t)$ approximated as $m(i\Delta x,k\Delta t)=m_{i,k}=n_{i,k}/\mathcal{N}$, with $\sum_{i=1}^{L/\Delta x}n_{i,1}=\mathcal{N}$; MATLAB functions for the GRW solvers that are called in the main code, describing the evolution of the groups of particles $n_{i,k}$; auxiliary MATLAB functions, e.g., to plot the results or to compute analytical solutions and various source terms occurring in the equations. This formulation facilitates the comparison with other discretization schemes presented in this article. However, when we are interested in describing the evolution of small numbers of ingressing particles, the same codes can be used to solve for $n_{i,k}$, without the need to convert the number of particles to approximations $m_{i,k}$ of the continuous field $m(x,t)$. %\red{$\leftarrow\leftarrow$}%\red{new text: $\rightarrow\rightarrow$}
[The GRW algorithms  presented in this article are  implemented in MATLAB and  are openly available in the Git repository
\href{https://github.com/PMFlow/Stefan_type_Problems} {https://github.com/PMFlow/Stefan\_type\_Problems}.]
%\red{$\leftarrow\leftarrow$}\red{(... to be completed soon; Should we include the FEM codes used in Sects. 3.1 and 4.3?)}

\section{Benchmarking the GRW method for handling moving interfaces}
\label{sec:GRWStefan}
\subsection{The case of the classical Stefan problem}
\label{sec:Stefan}
A Stefan problem is a specific initial boundary value problem describing the heat distribution in a phase-changing medium. It usually needs to be solved in a time-changing
space domain with  moving boundary conditions.  As the moving boundary is a function of time and it has to be determined as a part of the unknowns of the problem, the exact solution to free- and moving-boundary problems is rarely available. Therefore, various numerical methods have been developed to approximate the solution to the Stefan problem, for instance, finite difference method \cite{kutluay1997numerical, savovic2003finite}, finite element method \cite{mori1977finite, mori1976stability}, probabilistic approaches \cite{casaban2023numerical,Ogren2022}.
In this section, we examine the application of the  GRW method for solving the classical  Stefan problem.  We begin with providing a brief overview of the problem and its parameters. We then present the GRW simulation results and compare them against both the analytical solution (see for instance,  ~\cite{Gupta2003, Stefan1891}) and the random walk (RW) approximation from~\cite{Ogren2022}. Finally, we study the grid convergence behaviour in the GRW method.  
%and compare GRW with numerical RW results from~\cite{Ogren2022}, and to the famous analytic solution of the classical Stefan problem~\cite{Gupta2003, Stefan1891}.

\subsubsection{Problem formulation}
\label{sec:formStefan}

%GRW solutions\footnote{See \cite[Sect. 3.3]{Suciu2019} for the description of the GRW algorithms.} to the one phase Stefan's problem \cite[Eqs. (16-21)]{Ogren2022} are compared in the following to the analytical solution Eq.~(23) in the cited reference. \\

We start off  this section by introducing the classical Stefan problem and briefly explain the model's parameters.  
The classical Stefan problem describes the solidification and melting process of ice. Initially, the ice is solid, so there is no liquid phase. We consider the temperature to be zero at the solid phase. For $t > 0$ the ice can start to melt and thus we can have a
water phase, i.e. a liquid phase. Let  $m(x, t)$ denote the temperature distribution in the liquid phase at time $t$ and at the position $x$. Let $x = s(t)$ denote the position of the moving interface that separates the solid and liquid phases, while $m(x, t)$  is defined in the region $Q_s(T)$.
%defined by $$ Q_s(T):= \{ (x, t) \;   | \;    x \in (0, s(t))  \; \text{and}\;  t \in (0, T)\}.$$ 
The problem reads: Find the temperature profile
%the diffusant 
$m(x, t)$ and the position of the moving interface  $x = s(t)$ for $t\in(0, T)$ such that the couple $(m(x, t), s(t))$ satisfies the following set of equations:
\begin{align}
\label{stefan2a11}&\displaystyle \frac{\partial m}{\partial t} -\alpha_L \frac{\partial^2 m}{\partial x^2} = 0\;\;\; \ \text{in}\;\;\; Q_s(T),\\
\label{stefan2a12}
&  m(0, t) = m_L  \;\;\; \text{for}\;\;  t\in(0, T),\\
\label{stefan2a13}
&  m(s(t), t) =0  \;\;\; \text{for}\;\;  t\in(0, T),\\
\label{stefan2a14}&  -k_L \frac{\partial m}{\partial x}(s(t), t) = \rho_L \ell s^{\prime}(t) \;\;\; \text{for}\;\; t\in(0, T),\\
\label{stefan2a15}&m(x, 0) = 0, \;\; x\in [0,s(0)],\\
\label{stefan2a16} & s(0) = s_0\geq 0.
\end{align}
The parameter $\alpha_L$ represents the thermal diffusivity of the liquid and is defined by 
\[\alpha_L := \frac{k_L}{\rho_L c_L}. \] Here $k_L>0$ is the corresponding heat conductivity, $\rho_L>0$ is the density and $c_L>0$  is the specific heat capacity of the liquid phase.  Set $\ell>0$ as the specific latent heat required to melt ice. 
Assuming that the densities in the liquid and solid phases are equal,
the analytical solution to~\eqref{stefan2a11}--\eqref{stefan2a16} is given by~\cite{Stefan1891}
\[\widetilde{m}(x, t) = m_L \left( 1 - \frac{{\rm erf}(x/(2 \sqrt{\alpha_L t}))}{{\rm erf}(\lambda)}\right), \] and 
\[\widetilde{s}(t) = 2 \lambda \sqrt{\alpha_L t},\]
where $m_L>0$ is a given constant and $\lambda$ satisfies the following algebraic equation
\[\mu \sqrt{\pi} \lambda e^{\lambda^2} {\rm erf}(\lambda)= m_L,\] 
with $\mu := \ell/c_L$ and $s_0=0$.
We refer the reader to~\cite{Gupta2003, Stefan1891} for detailed descriptions of the Stefan problem and to~\cite{Ogren2022} for a random walk method approximating the solution of this problem.  In the following we solve problem~\eqref{stefan2a11}--\eqref{stefan2a16} in dimensionless form, with $\alpha_L=1$ and $m_L=1$, for the final time $T=0.5$.

\subsubsection{GRW convergence study for a benchmark case}

 The Stefan problem described in Section~\ref{sec:formStefan} is solved now with a GRW scheme based on the algorithm described in Section~\ref{sec:grw} using the parameters of the numerical example presented in \cite{Ogren2022}. 
 %Here and in the following, the GRW solutions are computed with the fixed number of particles $\mathcal{N}=10^{24}$, which ensures the convergence of the scheme.% (see Remark~\ref{rem:GRWconv} ??????????????????????????????????????).

\begin{figure}[h]
%\begin{minipage}[c]{0.48\linewidth} %\centering
\includegraphics[width=0.48\linewidth]{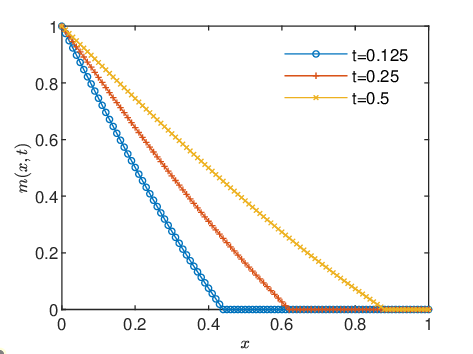}
%\caption{\label{fig:c_Stefan}{\footnotesize Water temperatures as functions of space for successive times.}}
%\end{minipage}
\hspace{0.1cm}
%\begin{minipage}[c]{0.48\linewidth} %\centering
%\includegraphics[width=0.48\linewidth]{eps/s_Stefan.eps}
\includegraphics[width=0.42\linewidth]{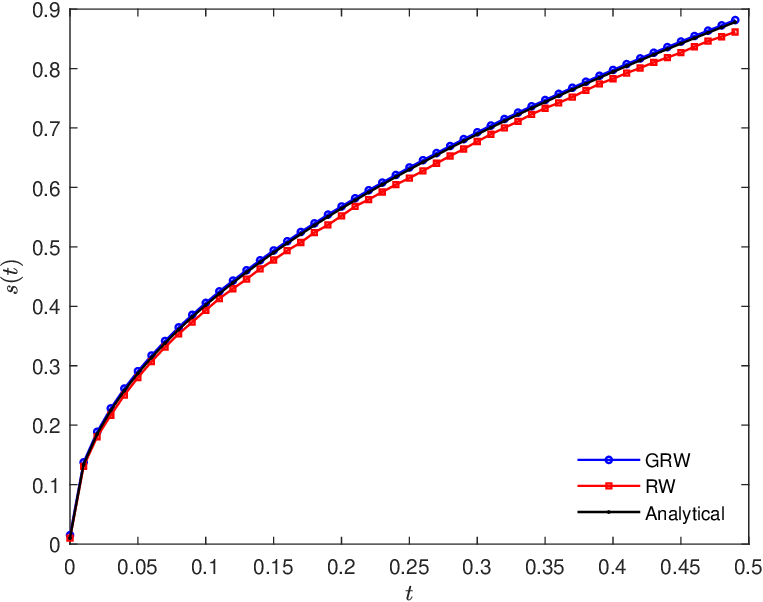}
\caption{\label{fig:s_Stefan}{{\em Left}: Profiles $m(x,t)$ approximated by GRW method. {\em Right}: Comparison of the GRW simulation results with the RW and the analytical solutions  for the position of the  moving boundary $s(t)$.}}
%\end{minipage}
\end{figure}

The left-hand side of Figure~\ref{fig:s_Stefan} shows the distribution of the water temperature at successive times, and in the right-hand side of Figure~\ref{fig:s_Stefan} the GRW solution for the moving interface water/ice is compared with the exact analytical solution and the solution obtained with a random walk (RW) method presented in \cite[Fig. 8]{Ogren2022}. The two numerical solutions were computed with the same spatial discretization $\Delta x = 0.01$. The relative errors of the RW and GRW solutions with respect to the analytical solution, computed in terms of $\ell^2$ vector norms, were of $1.86 \%$ and $0.56 \%$, respectively.

%\red{(see /Matlab/1D/Stefan/main\_Stefan\_conv.m and .../plot\_s.m)}

Following the common approach of grid-convergence tests cf. e.g. \cite{Suciuetal2021,Suciuetal2024,Roy2005,Alecsaetal2020}, we assess the convergence of the numerical solution $(m(x,t),s(t))$ to the analytical solution $(\widetilde{m}(x,t),\widetilde{s}(t))$ in terms of  estimated orders of convergence ($EOC$) given by the slope of the decay of the errors in the logarithmic scale after successively halving the discretization parameter $\Delta x$, i.e.
\begin{equation}\label{eq:eoc}
EOC:=\log\left(\frac{e_l}{e_{l+1}}\right)/\log(2),\;\; l=1,2,\cdots,
\end{equation}
where $e_\ell$ is the error measured in discrete  $L^2$ norm~\cite{Alecsaetal2020}  %%(i.e., $\|s_l-\widetilde{s}_l\|_{L^2(0, T)}$ ) 
(i.e. in terms of  $\| \cdot \| := \sqrt{\Delta x} \| \cdot \|_{\ell^2}$, where $\| \cdot \|_{\ell^2}$ is the Euclidean vector norm, in $[0,s(T)]$, for $m(x,t)$, and $\| \cdot \|:=\sqrt{\Delta t} \| \cdot \|_{\ell^2}$ in $[0,T]$, for $s(t)$).
%for the position of moving boundary and measured in discrete $L^2(0, T, L^2(0, s(t)))$ norm (i.e., $\|m_l-\widetilde{m}_l\|_{L^2(0, T, L^2(0, s(t)))}$) for the profile  \cite{Alecsaetal2020}. 
%\red{I do not understand the norm notations: in the code the error norms of the ``profile'' $m(x,t)$ are computed in [0,1]; both the numerical and the analytical solutions have the length of the vector $x=0:dx:1$; both solutions vanish for $x_i>s(t)$. Then, effectively the norm is computen in $[0,s(t)]$, with $t\in [0,T]$.} 
The tables presenting the convergence tests in the following have a generic structure:
a column for the spatial discretization parameter $\Delta x$, columns with errors of the profile $m(x,t)$ and of the position of the moving boundary $s(t)$ and the corresponding EOC, and another column with the computing times (CT) for successive $\Delta x$.

The results of the convergence test and the corresponding computing times (CT) for the GRW approximation of the solution to the Stefan problem are presented in Table~\ref{tab:conv_c_s}. % \red{(see /Matlab/1D/conv\_Stefan/)}

\begin{table}[h]
\begin{center}
\begin{tabular}{ c c c c c c c c}
  \hline
   & $\Delta x$ & $\|m-\widetilde{m}\|$ & $EOC$ & $\|s-\widetilde{s}\|$ & $EOC$ &  CT (sec) &\\
  \hline
   &  4.00e-02  &  6.52e-03 & --   &  9.62e-03 & --   &  0.04 &\\
   &  2.00e-02  &  3.49e-03 & 0.90 &  4.89e-03 & 0.98 &  0.05 &\\
   &  1.00e-02  &  1.84e-03 & 0.92 &  2.41e-03 & 1.02 &  0.08 &\\
   &  5.00e-03  &  9.22e-04 & 1.00 &  1.18e-03 & 1.03 &  0.32 &\\
   &  2.50e-03  &  4.28e-04 & 1.11 &  5.76e-04 & 1.03 &  1.84 &\\
  \hline
\end{tabular}
\end{center}
\caption{\label{tab:conv_c_s} %\red{new - } 
%\red{(The table caption, here and in the following, should be reformulated; ``mesh size'',  ``errors'', and the acronym EOC have been introduced in Eq. (7).)} 
%Mesh size, errors, estimated orders of convergence (EOC), and computing time (CT) of the GRW approximation for the classic Stefan problem.
Orders of convergence and computing times of the GRW approximation for the 
%classic 
classical Stefan problem.}
\end{table}

%\begin{table}[h]
%\begin{center}
%\begin{tabular}{ c c c c c c c c}
%  \hline
%   & $\Delta x$ & $\|c-\widetilde{c}\|$ & $EOC$ & $\|s-\widetilde{s}\|$ & $EOC$ &  CT (sec) &\\
%  \hline
%   &  1.00e-02  &  1.84e-03 & --   &  2.41e-03 & --   &  0.13 &\\
%   &  5.00e-03  &  9.22e-04 & 1.00 &  1.18e-03 & 1.03 &  0.36 &\\
%   &  2.50e-02  &  4.28e-04 & 1.11 &  5.76e-04 & 1.03 &  2.04 &\\
%   &  1.25e-03  &  2.15e-04 & 0.99 &  2.82e-04 & 1.03 & 13.60 &\\
%   &  6.25e-04  &  1.00e-04 & 1.10 &  1.39e-04 & 1.02 & 94.43 &\\
%  \hline
%\end{tabular}
%\end{center}
%\caption{\label{tab:conv_c_s} %\red{new - } 
%Orders of convergence of the GRW solution $(c(t, x), s(t))$ and %computational times (CT).}
%\end{table}

In the following, we compare the GRW approximation of the solution to the Stefan problem with approximations (of the same solution) obtained using the finite element method (FEM).
To proceed with the FEM approach, we first transform the problem formulated in the moving domain $(0, s(t))$ into one posed for a fixed domain $(0,1)$. To this end, we employ the Landau transformation $y = x/s(t)$. We then solve the transformed equations discretizing the fixed domain $(0, 1)$ with a uniform grid size $\Delta y$. We refer the reader to \cite{nepal2021moving} for more details on how FEM works for one-dimensional moving-boundary problems. 

%A comparison of the analytical and FEM solutions is shown in Figs.~\ref{fig:stefan_tem_fem}~and~\ref{fig:stefan_mb_fem}. 
The estimated orders of convergence and the corresponding computing times are presented in Table~\ref{tab:conv_ana_FEM}. %\green{At this time, we calculate the error $e_{l}$ using the discrete $L^2(0, T; L^2(0, 1))$  norm for the concentration profile $m(y, t)$, as defined in \cite{nepal2021error}  and the discrete $L^2(0, T)$ norm, as defined above, for the position of the moving boundary $s(t)$}. 
Comparing with Table~\ref{tab:conv_c_s} we see that both GRW and FEM solutions converge with the same order $EOC\approx 1$. However, the GRW method, which solves the moving-boundary problem directly, without transformation to a fixed domain, approximates the solution with an accuracy that is systematically higher than that of the FEM solution and the corresponding computing times are one or two orders of magnitude smaller. 

%\begin{figure}[h]
%\begin{minipage}[c]{0.48\linewidth} %\centering
%\includegraphics[width=\linewidth]{stefan_problem_tem.png}
%\caption{\label{fig:stefan_tem_fem}{\footnotesize Water temperatures %as functions of space at time $t=0.5$.}}
%\end{minipage}
%\hspace{0.1cm}
%\begin{minipage}[c]{0.48\linewidth} %\centering
%\includegraphics[width=\linewidth]{stefan_problem_mb.png}
%\caption{\label{fig:stefan_mb_fem}{\footnotesize The moving boundary %as a function of time.}}
%\end{minipage}
%\end{figure}

\begin{table}[h]
\begin{center}
\begin{tabular}{ c c c c c c c c}
  \hline
   & 
   %$\Delta y$
   $\Delta y$ 
   & $\|m-\widetilde{m}\|$ & $EOC$ & $\|s-\widetilde{s}\|$ & $EOC$ & CT (sec) &  \\
  \hline
   &  4.00e-02  & 2.27e-02 &   --  &2.32e-02  &  --   & 0.41 & \\ %0.25 & \\
   &  2.00e-02  & 1.11e-02 & 1.03  & 1.16e-02 & 1.00 &  1.06 &  \\ %0.43 &  \\
   &  1.00e-02  & 5.49e-03 & 1.02  & 5.77e-03 & 1.01 &  2.46 & \\ %2.90 & \\
   &  5.00e-03  & 2.73e-03 & 1.01  & 2.85e-03 & 1.02 & 20.73 & \\ %16.40  & \\
   &  2.50e-03  & 1.36e-03 &  1.00 & 1.40e-03 & 1.03 &  203.07 & \\ % 128.49 & \\
  \hline
\end{tabular}
\end{center}
\caption{\label{tab:conv_ana_FEM}%\red{Please complete for the convergence of c(t) by halving 4 times the initial dx=1e-2 !} 
% Mesh size, errors, estimated orders of convergence (EOC), and computing time (CT) of the FEM approximation for the classic Stefan problem.
Orders of convergence and computing times of the FEM approximation for the 
%classic 
classical Stefan problem.}
\end{table}

\subsection{The case of the Stefan problem with kinetic condition}
\label{sec:Stefan_Kinetic}

Now, we explore the convergence behaviour of GRW for the case of a Stefan problem with a kinetic condition. To begin with, we introduce a prototype model. Then, we construct an analytical solution which is further used in code verification and convergence tests.  
%%\subsection{Problem formulation}
%\subsection{Convergence study for a constructed solution}
%%\section{Constructed solution (Stefan-type problem)}
\subsubsection{Problem formulation}
\label{Stefan_Kinetic_form}

Consider that $m(x, t)$  acts in the region $Q_s(T)$. 
%defined by $$ Q_s(T):= \{ (x, t) \; | \; x \in (0, s(t))  \; \text{and}\; t \in (0, T)\}.$$ 
The problem reads: Find the profile 
%the diffusant 
$m(x, t)$ and the position of the moving interface $x = s(t)$ for $t\in(0, T)$ such that the couple $(m(x, t), s(t))$ satisfies the following:

\begin{align}
\label{2a11}&\displaystyle \frac{\partial m}{\partial t} -D \frac{\partial^2 m}{\partial x^2} = 0\;\;\; \ \text{in}\;\;\; Q_s(T),\\
\label{2a12}
&  m(0, t) = m_0(t)  \;\;\; \text{for}\;\;  t\in(0, T),\\
\label{2a13}&-D \frac{\partial m}{\partial x}(s(t), t)  = s^{\prime}(t)\,m(s(t), t)  \;\;\; \text{for}\;\; t\in(0, T),\\
\label{2a14}&s^{\prime}(t) = D\alpha\left(m(s(t), t)-\alpha s(t)\right) \;\;\;\;\text{for } \;\;\; t \in (0,T),\\
\label{2a15}&m(x, 0) = 0 \;\;\;\text{for}\;\;\; x \in [0, s(0)],\\
\label{2a16} & s(0) = s_0>0.
\end{align}

We notice that, unlike in the classical Stefan problem formulated in Section~\ref{sec:formStefan}, the interface condition (\ref{2a13}) is now 
%constraint 
constrained 
by the kinetic condition (\ref{2a14}) which also renders the coupled problem nonlinear.

%For the choices of $\alpha = 0.1, s_0 = 0.001$ and $D = 1$, the plots %for finite element approximation and analytical solutions, $m(t, x)$ %at time $t= 10$ and $s(t)$ is shown in Figure \ref{fig:m_Stefan_type} %and Figure \ref{fig:st_Stefan_type} respectively. 
%\begin{figure}[h]
%\begin{minipage}[c]{0.48\linewidth} %\centering
%\includegraphics[width=\linewidth]{m_stefan_type.png}
%\caption{\label{fig:m_Stefan_type}{\footnotesize $m(t, x)$ at $t=10$.}}
%\end{minipage}
%\hspace{0.1cm}
%\begin{minipage}[c]{0.48\linewidth} %\centering
%\includegraphics[width=\linewidth]{st_stefan_type.png}
%\caption{\label{fig:st_Stefan_type}{\footnotesize Position of the %moving boundary.}}
%\end{minipage}
%\end{figure}
%\AM{ Surendra: How does it compare against your FEM approximation? Can you please check this and show us a pointwise comparison? At a first sight, for your manufactured example, $T_f$ must be finite and it cannot be too big as it must preserve $s(t)\geq 0$. }
%----------------------------------------

\subsubsection{Convergence study for a benchmark case}
\label{Stefan_Kinetic_conv}

In the absence of analytical solutions, the code verification and the grid-convergence test can be done by computing errors with respect to manufactured solutions. The latter are functions $\widetilde{u}$ of the dependent variables which have to fulfill the only requirement that the result of applying the differential operator $P$, of the partial differential equation $Pu=0$, is a non-vanishing function $f$. Then, $\widetilde{u}$ is an exact solution of the original partial differential equation with the source term $f$ added in the right-hand side, i.e., $P\widetilde{u}=f$~\cite{Roache2002}. 

In the following, we choose the analytical solution 
\begin{align}
\widetilde{m}(x, t) &= \sqrt{\frac{t}{T}}\,\exp(-\alpha \, x), \  \alpha > 0, \nonumber\\
\widetilde{s}(t) &= s_0 + 2 D \alpha \sqrt{\frac{t}{T}}. \label{eq:man_sol_St_kin}
\end{align}
Inserting \eqref{eq:man_sol_St_kin} into the model \eqref{2a11}--\eqref{2a16} one obtains the source terms
\begin{align*}
f_1(x,t)&=\frac{1}{2\sqrt{t\,T}}\exp(-\alpha x) - D \alpha^2 \sqrt{\frac{t}{T}} \exp(-\alpha x),\\
f_2(\widetilde{s}(t),t)&=(D\alpha -\widetilde{s}^{\,\prime}(t))\widetilde{m}(\widetilde{s}(t),t),\\
f_3(t)&=D\alpha\left(\frac{1}{\sqrt{t\,T}}-\widetilde{m}(\widetilde{s}(t),t)+\alpha \widetilde{s}(t)\right),
\end{align*}
in the right-hand sides of
%Eqs.
~(\ref{2a11}), (\ref{2a13}), and (\ref{2a14}), respectively. Also, the manufactured solution (\ref{eq:man_sol_St_kin}) fixes the boundary condition (\ref{2a12}) to $m_0(t)=\sqrt{t/T}$ and verifies the initial conditions (\ref{2a15}) and (\ref{2a16}).

%As shown by Eqs.~(\ref{2a13})-(\ref{2a14}), in this case the problem is nonlinear. %\red{new text: $\rightarrow\rightarrow$}

\begin{table}[h]
\begin{center}
\begin{tabular}{ c c c c c c c c }
  \hline
   & $\Delta x$ & $\|m-\widetilde{m}\|$ & $EOC$ & $\|s-\widetilde{s}\|$ & $EOC$ &  CT (sec) & \\
  \hline
   &  1.00e-02  & 1.39e-03 & --   & 1.03e-03 & --   & 0.33 & \\
   &  5.00e-03  & 4.99e-04 & 1.48 & 5.15e-04 & 1.00 & 1.12 & \\
   &  2.50e-03  & 1.78e-04 & 1.49 & 2.58e-04 & 1.00 & 7.13 & \\
   &  1.25e-03  & 6.32e-05 & 1.49 & 1.29e-04 & 1.00 & 44.91 & \\
   &  6.25e-04  & 2.24e-05 & 1.50 & 6.45e-05 & 1.00 & 376.95 & \\
  \hline
\end{tabular}
\end{center}
\caption{\label{tab:conv_Stefan_kin}
%%Mesh size, errors, and estimated orders of convergence (EOC)  of the GRW approximation for the Stefan problem with the kinetic condition.
Orders of convergence and computing times of the GRW approximation for the Stefan problem with the kinetic condition.}
\end{table}

Preliminary tests with a linearization $L$-scheme \cite{Suciuetal2021} have shown that the nonlinearity $m^2(s(t),t)$ arising when (\ref{2a14}) is inserted in (\ref{2a13}) is solved after the first iteration. Therefore, in solving Stefan problems with a kinetic condition we use the so-called explicit linearization of the moving boundary condition, consisting of replacing $m(s(t),t)$ by $m(s(t-\Delta t),t)$ in the right-hand side of (\ref{2a13}). When, as in the present case, we evaluate the convergence with respect to an exact solution, the moving boundary condition (\ref{2a13}) is also formulated exactly by using $\widetilde{m}(\widetilde{s}(t),t)$. [For details on the implementation see the code 'main\_Stefan\_kin\_conv.m' in the folder 'Stefan\_kinetic' of the Git repository.]%\red{$\leftarrow\leftarrow$}. 

The results obtained with $\alpha = 0.1, \; D = 1, \; s_0 = 0.1$, and $T=1$ are presented in Table~\ref{tab:conv_Stefan_kin}. %\red{(see /Matlab/1D/conv\_Stefan\_kin/)}\\

 Observe that, when approximating $m(x,t)$, the convergence orders are improved  compared to what is normally achieved for the classical Stefan problem. 

%{\bf what does the reader see in the Table? please explain}

%\red{- results for T=10 and T=1}\\
%\red{- should we present results for T=1, with EOCc ~ 1.5 and EOCs = %0.997, 0.998, 0.999, instead of those for T=10 with the huge %EOCc=12.4?}
%\blue{In my opinion, it is okay to present results for T=1.}\\
%\red{04.07.2024: In Table 3 are now the results for T=1 !}

\section{GRW approximation capturing diffusion in rubber}
\label{sec:rubber}

Within this section we present GRW simulation results for a moving-boundary problem describing the penetration of diffusants into rubbers. Our interest in this particular scenario is linked with the many potential technological applications (cf. e.g. \cite{nepal2021moving, lewis1980laboratory, hayes1955diffusion,  kaliyathan2021analysis}), where we believe we can contribute with numerical predications of the durability (service life) of the rubber-based material.   

We start by introducing and explaining the model equations. As a next step, we study the convergence behaviour for a benchmark case where we perform the grid-convergence test for a problem with the analytical solution. Finally, we present the GRW simulation results and compare them with the FEM and random walk (RW) solution for our model equations. 
\subsection{Problem formulation}
\label{sec:rubber_form}

We recall the following problem from \cite{nepal2021moving}: For a fixed observation time $T\in (0, \infty)$, the interval $[0, T]$ represents the duration of the involved physical processes. Let $t \in [0, T]$ denote the time variable, $s(t)$ the position of the moving-boundary at time $t$ and  $x\in [0, s(t)]$ the space variable. The function $m(x, t)$ represents the concentration of diffusant placed at position $x$ and time  $t$.
The diffusants concentration  $m(x, t)$  acts in the 
%non-cylindrical parabolic 
region $Q_s(T)$.
%defined by $$ Q_s(T):= \{ (x, t) \; | \;  x \in (0, s(t)) \; \text{and}\; t \in (0, T) \}.$$ 
The problem is: Find 
the diffusant concentration profile  
$m(x, t)$  and  the position of the moving-boundary $x = s(t)$ such that the couple $(m(x, t), s(t))$ satisfies the following:
\begin{align}
& \frac{\partial m}{\partial t} -D \frac{\partial^2 m}{\partial x^2} = 0 \;\;\; \ \text{in}\;\;\; Q_s(T),\label{Eq:r1}\\
& -D \frac{\partial m}{\partial x}(0, t)  = \beta(b-{\rm H}m(0, t)) \;\;\; \text{for}\;\; t\in(0, T),\label{Eq:r2}\\
& -D \frac{\partial m}{\partial x}(s(t), t)  = s'(t)m(s(t), t) \;\;\; \text{for}\;\; t\in(0, T),\label{Eq:r3}\\
& s^{\prime}(t)=a_0\left(m(s(t), t)-  \sigma(s(t))\right) \;\;\; \text{for}\;\; t\in(0, T),\label{Eq:r4}\\
& m(x, 0) = m_0(x)\;\;\;\text{for}\;\;\; x \in [0, s(0)],\label{Eq:r5}\\
& s(0)=s_0>0.\label{Eq:r6}
\end{align}
Here $D$ is the diffusion coefficient, $a_0>0$ is a kinetic coefficient, $\beta$ is a positive constant  
and  $\text{H}>0$ is called the Henry constant. Additionally, $\sigma$ is a function on $\mathbb{R}$ describing the swelling of the rubber, $b$ represents the amount of diffusant concentration available at the boundary.  We consider the diffusion equation in \eqref{Eq:r1} to represent the role of diffusion for the evolution of the diffusant concentration $m(x, t)$. At the left boundary $x=0$,  we define the Robin-type boundary condition \eqref{Eq:r2} to describe the inflow of diffusants into rubber.  At the right boundary $x = s(t)$, we have the Robin boundary condition \eqref{Eq:r3} that describes the mass conservation of the diffusants at the moving boundary. The ordinary differential equation  \eqref{Eq:r4} describes the speed of the moving boundary. In  \eqref{Eq:r6}, $s_0$ denotes the initial position of the moving-boundary, while in \eqref{Eq:r5} $m_0$ denotes the concentration of the diffusant at $t=0$.

The model equations \eqref{Eq:r1}--\eqref{Eq:r6} were designed to mimic experimental results involving
the ethylene propylene diene monomer rubber being swollen in
cyclohexane. We do not discuss here why this particular choice of materials is relevant and prefer instead to indicate a couple of mathematical results that delimitate the mathematical validity of the model. We refer the reader to Theorem 3.4 in \cite{kumazaki2020global} for detailed results on the global existence of weak solutions to \eqref{Eq:r1}--\eqref{Eq:r6} and continuous dependence with respect to initial data. The detailed description of model equations, laboratory experimental setup, and finite element simulation results is presented in \cite{nepal2021moving}. The finite element approximation recovered the experimental data well. To approximate the solution using the FEM, we first transform the moving domain $(0, s(t))$ into a fixed domain $(0,1)$ by using the transformation $y = x/s(t)$. We then solve the transformed equations by discretizing the fixed domain $(0, 1)$ with a uniform grid size $\Delta y$.  The convergence analysis of the finite element scheme is provided in \cite{nepal2023analysis,nepal2021error} where we achieved first-order convergence in both space and time for the position of the moving boundary and diffusant profile. 

\subsection{Convergence study for a benchmark case for diffusion in rubber}
\label{sec:rubber_conv}
%%The analytical solution to \eqref{Eq:r1}-\eqref{Eq:r6} is unavailable.

In the absence of analytical solutions, the code verification and the grid-convergence test can also here be done by computing errors with respect to manufactured solutions.
%
%%Instead of solving \eqref{Eq:r1}-\eqref{Eq:r6}, we assume a solution that does not satisfy the governing equations  \eqref{Eq:r1}-\eqref{Eq:r6}. This assumed solution is then inserted into \eqref{Eq:r1}-\eqref{Eq:r6}, resulting in the residual term, which is added back to the equations \eqref{Eq:r1}-\eqref{Eq:r6} as a source term. The modified problem then has the assumed solution as an exact solution and errors can be computed based on the exact solution, i.e., manufacture solution to confirm the order of accuracy of the GRW method.  
%
%\red{(Surendra, if you constructed a manufactured solution for the RWM - rubber problem, please put it in the Dropbox, so that I can use it to assess the convergence of the GRW solution. We should compare the two RW approaches.)}
From now onwards,  we choose a linear dependence of  $\sigma$ on $s(t)$ in \eqref{Eq:r4}.   % Eq.~(\ref{Eq:r4}). 
With this choice, \eqref{Eq:r4} becomes 
\begin{equation}\label{Eq:r4_lin}
 s^{\prime}(t)=a_0\left(m(s(t), t)-\varrho s(t)\right) \;\;\; \text{for}\;\; t\in(0, T),
\end{equation}
where $\varrho$ is a constant parameter. This linear choice of $\sigma$ is motivated by our previous work \cite{nepal2021moving,Nepaletal2023} where we explored the numerical simulation results for three different choices of $\varrho$ and $a_0$. 
We choose the following manufactured solution for the concentration profile and the diffusion front:
\begin{align}
\widetilde{m}(x, t)& = \left(1-\frac{x}{\widetilde{s}(t)}\right)^3\cos\left(\frac{t}{T}\right),\nonumber\\
\widetilde{s}(t)& = s_0\left(\frac{t}{10T}-\frac{t^2}{10T^2}+\frac{t^3}{30T^3}+1\right).\label{eq:man_sol}
\end{align}
The exact solution \eqref{eq:man_sol} solves a modified problem, with the following source terms added in the right-hand side of 
%Eqs.
~\eqref{Eq:r1}--\eqref{Eq:r4}:
%\begin{align*}
%f_1(x, t) = &
%\frac{\sin\left(t/T\right)}{T}  \left( \frac{x}{\widetilde{s}(t)}  -  1 \right)^3  + 
%\frac{3\,x\,\cos\left(t/T\right)}{\left( \widetilde{s}(t) \right)^2}   s_m^{\prime}(t)
%\left( 1 - \frac{x}{\widetilde{s}(t)} \right)^2 + 
% \frac{6\,D\,\cos\left(t/T\right)}{\left( \widetilde{s}(t)\right)^2} \left( \frac{x}{\widetilde{s}(t)}  -  1 \right), \\
%f_2(t) =  &-b \beta +\beta  H \cos \left(\frac{t}{T}\right) + \frac{3 D \cos(t/T)}{\widetilde{s}(t)} , \\
%f_3(t) = & \ 0 , \\
%f_4(t) = & \widetilde{s}^{\prime}(t) + 
%\frac{a_0}{\varrho} \widetilde{s}(t).
%\end{align*}
%\red{Dear Yosief, your interpretation of the source terms is more compact and understandable. I only replaced the (old) $s_{m}^{\prime}$ by $\widetilde{s}^{\prime}$, rearranged some terms, and I inserted some explanations. Please check and keep one version.}
% (YW): I think we keep this formulation of yours (NS).
\begin{align*}
f_1(x, t) = & 
\ \frac{\sin\left(t/T\right)}{T}  \left( \frac{x}{\widetilde{s}(t)}  -  1 \right)^3  + 
\frac{3\,x\,\cos\left(t/T\right)}{\left( \widetilde{s}(t) \right)^2}  
\left( \frac{x}{\widetilde{s}(t)} -1\right)^2\widetilde{s}^{\,\prime}(t) + 
 \frac{6\,D\,\cos\left(t/T\right)}{\left( \widetilde{s}(t)\right)^2} \left( \frac{x}{\widetilde{s}(t)}  -  1 \right), \\
f_2(0,t) =  & \ -\beta\, b +\left(\beta \,  H  + \frac{3 D }{\widetilde{s}(t)}\right)\cos\left(\frac{t}{T}\right) , \\
f_3\left(\widetilde{s}(t),t\right) = & \ \ 0 , \\
f_4(t) = & \ \widetilde{s}^{\,\prime}(t) + 
\frac{a_0}{\varrho} \widetilde{s}(t).
\end{align*}
While deriving the above source terms, we made use of the following properties of the manufactured solutions:
\begin{align*}
&\widetilde{m}(0, t)=\cos\left(\frac{t}{T}\right),\; \frac{\partial\widetilde{m}}{\partial x}(0, t)=\frac{3}{\widetilde{s}(t)}\cos\left(\frac{t}{T}\right), \\
&\widetilde{m}\left(\widetilde{s}(t), t\right)=0,\; \frac{\partial\widetilde{m}}{\partial x}\left(\widetilde{s}(t), t\right)=0.
\end{align*}
The initial conditions of the modified problem, corresponding to (\ref{Eq:r5}) and (\ref{Eq:r6}), are specified by the manufactured solution (\ref{eq:man_sol}): $m_0(x)=\widetilde{m}(x, 0)$, $s_0=\widetilde{s}(0)$.

%\red{(see code /Matlab/1D/conv\_rubber/main\_RobinBC\_conv.m)}\\
%\red{- corected sources f1-f4, with s(t) instead of mx, which was the %constant value of s at fixed t, not a function!}\\
%\red{- computations without function\_handle Ec = $@...$, which breaks %the convergence of c}\\
%\red{- results for Dirichlet and Robin BC}

%\newpage  % for the time being

We approach this problem (in dimensionless form) numerically via the GRW method. We make use   of the following parameters: $D=1,\; \beta=1,\; b=1,\; H=1,\; a_0=1,\; \varrho=1,\; s_0=0.1$, and T=0.001. %\red{new text: $\rightarrow\rightarrow$} 
The condition on the moving boundary is formulated exactly with the aid of the manufactured solution, as in Section~ \ref{Stefan_Kinetic_conv} above. %\red{$\leftarrow\leftarrow$}
We compute the error of the concentration profile and moving front by comparing the GRW approximation with the manufactured solution in terms of the discrete $L^2$ norms introduced below equation~(\ref{eq:eoc}). 
%%$L^2(0, T, L^2(0, s(t)))$ 
%%%%\red{Please check: should it be $L^2(0, s(t))$?}
%%and $L^2(0, T)$ norms,  respectively.  
The obtained errors and estimated orders of convergence are shown in Table~\ref{tab:conv_rubber}.

%\begin{table}[h]
%\begin{center}
%\begin{tabular}{ c c c c c c c }
%  \hline
%   & $\Delta x$ & $\|c-c_m\|$ & $EOC$ & $\|s-s_m\|$ & $EOC$ &  \\
%  \hline
%   &  1.00e-02  & 1.78e-03 & --   & 5.82e-05 & --   & \\
%  &  5.00e-03  & 4.21e-04 & 2.08 & 1.44e-05 & 2.01 & \\
%   &  2.50e-03  & 1.03e-04 & 2.04 & 3.61e-06 & 2.00 & \\
%   &  1.25e-03  & 2.53e-05 & 2.02 & 9.01e-07 & 2.00 & \\
%   &  6.25e-03  & 6.30e-06 & 2.01 & 2.25e-07 & 2.00 & \\
%  \hline
%\end{tabular}
%\end{center}
%\caption{\label{tab:conv_rubber}Errors and  orders of convergence of the numerical GRW solution for the concentration profile $c(t, x)$ and the position of the moving boundary $s(t)$.}
%\end{table}

\begin{table}[h]
\begin{center}
\begin{tabular}{ c c c c c c c c }
  \hline
   & $\Delta x$ & $\|m-\widetilde{m}\|$ & $EOC$ & $\|s-\widetilde{s}\|$ & $EOC$ & CT (sec) & \\
  \hline
   &  1.00e-02  & 4.19e-03 & --   & 5.66e-06 & --   & 0.05 & \\
   &  5.00e-03  & 9.58e-04 & 2.13 & 1.43e-06 & 1.99 & 0.06 & \\
   &  2.50e-03  & 2.34e-04 & 2.03 & 3.58e-07 & 1.99 & 0.07 & \\
   &  1.25e-03  & 5.77e-05 & 2.02 & 8.99e-08 & 2.00 & 0.18 & \\
   &  6.25e-03  & 1.44e-05 & 2.01 & 2.25e-08 & 2.00 & 0.81 & \\
  \hline
\end{tabular}
\end{center}
\caption{\label{tab:conv_rubber}%%Mesh size, errors, and  estimated orders of convergence (EOC)  of the GRW approximation for the diffusant profile and for the position of the moving boundary.
Orders of convergence and computing times of the GRW approximation for the problem modeling diffusion in rubber.}
\end{table}
%\begin{remark}\label{rem:convOrders}
The performed grid-convergence tests indicate that the GRW approximations to the three Stefan-type problems converge with different rates: linear convergence for both concentration and penetration front in case of the classical Stefan problem (Table~\ref{tab:conv_c_s}), convergence of about an order of $1.5$ for concentration and linear convergence for the penetration front for the Stefan problem with kinetic condition (Table~\ref{tab:conv_Stefan_kin}), and second-order convergence of both concentration and penetration front for the model of diffusion in rubbers  (Table~\ref{tab:conv_rubber}). Having in view the rigorous numerical analysis results behind the proofs of convergence rates of FEM approximations for this class of problems (cf e.g. \cite{nepal2023analysis,nepal2021error} and references cited therein), we believe that the linear convergence for both concentration and penetration front observed in the case of the classical Stefan problem is what is in fact to be expected. The last two study cases indicate that, perhaps depending on the solution's regularity and choice of parameters, some sort of super-convergence may take place. It would be interesting to investigate  if one can delimitate rigorously parameter regimes and solution regularity classes where such super-convergence behavior can be guaranteed to happen.

\subsection{A forecast exercise}

%\red{What about the convergence in time? (see commented text in main.tex). Should we comment on it? Should we put the codes in the Git repository?}
%\red{For diffusion in rubber GRW convrges with EOC=2 in space and with EOC=1 in time, for both c(x,t) and s(t) (see code /Dropbox/drianMuntean/Matlab/1D/conv\_rubber/main\_RobinBC\_conv\_time.m).}

%\red{For Stefan problem GRW converges with EOC=1 in space and EOC=0.5 in time, for both c(x,t) and s(t) (see code /Dropbox/AdrianM-untean/Matlab/1D/conv\_Stefan/main\_Stefan\_conv\_time.m).} 

%\red{The difference between the convergence orders in space and time is similar to forward-time centered-space finite difference scheme for the parabolic advection-diffusion equation \cite{Strikwerda2004}, where the converge is second-order in space and first-order in time. We should try to understand why the convergence is slower for the the Stefan problem. Are there some analysis results on the convergence of Stefan problem?}

%\red{For Stefan problem with kinetic condition GRW converges with EOC$_c$=1.5 and EOC$_s$=1, in space, and similar EOC for the convergence in time (see code /Dropbox/AdrianMuntean/Matlab/1D/co-nv\_Stefan\_kin/main\_Stef\_kin\_conv\_time.m). This behavior should be investigated as well!}

After exploring the convergence of the GRW method for the problem with the manufactured solution in Section \ref{sec:rubber_conv}, our current goal is to validate the GRW method for \eqref{Eq:r1}--\eqref{Eq:r6}.  
As the analytical solution is unavailable, we compare the GRW simulation results with the FEM approximation, the RW solution, and laboratory experimental data presented in \cite{Nepaletal2023}. %\red{new text: $\rightarrow\rightarrow$} 
The moving boundary condition in the GRW code is now formulated with the aid of the explicit linearization described in Section \ref{Stefan_Kinetic_conv} above. [For details, see code `main\_RobinBC\_conv.m' in the folder `Rubber' of the Git repository).] %\red{$\leftarrow\leftarrow$}

The typical values for the chosen parameters and their dimensions used in simulations are provided in Table~\ref{GRWparameter}, based on \cite{Nepaletal2023}.

\begin {table}[ht]
\begin{center}
%	\resizebox{\columnwidth}{!}{%
 	\begin{tabular}{ |p{7.5cm}|p{2.0cm}|p{3.5cm}| }
%	\begin{tabular}{ |p{6.5cm}|p{2.0cm}|p{3.8cm}| }
		\hline
		Parameters & Dimension& Typical values\\
		\hline
		Diffusion constant, $D$ &$L^2T^{-1}$& $0.01$ (mm$^2$/min) \\\hline
		Absorption rate,  $\beta$ &$ LT^{-1}$ & $0.564$ (mm/min)\\\hline
		Constant,  $a_0$ & $L^4 T^{-1} M^{-1}$& $50$ (mm$^4$/min/gram) \\\hline
		Initial penetration front,  $s_0$ & $L$ & $0.01$ (mm)\\\hline
		$\sigma(s(t))$& $ML^{-3}$& $0.5s(t)$ (gram/mm$^3$)\\\hline
		Initial diffusant concentration, $m_0$ & $ML^{-3}$& $0.5$ (gram/mm$^3$)\\\hline
		Concentration in lower surface of the rubber, $b$& $ML^{-3}$&$10$ (gram/mm$^3$) \\\hline
		Henry's constant, 	$\rm{H}$& --& 2.50 (dimensionless)\\
		\hline
	\end{tabular}
%    }
	\caption {Parameters, their dimension, and choices of typical values for the reference model parameters.}
	\label{GRWparameter} 
\end{center}
\end {table}

\begin{figure}[ht]
%\begin{minipage}[c]{0.48\linewidth} %\centering
%\includegraphics[width=0.48\linewidth]{eps/c_rubber.eps}
\includegraphics[width=0.48\linewidth]{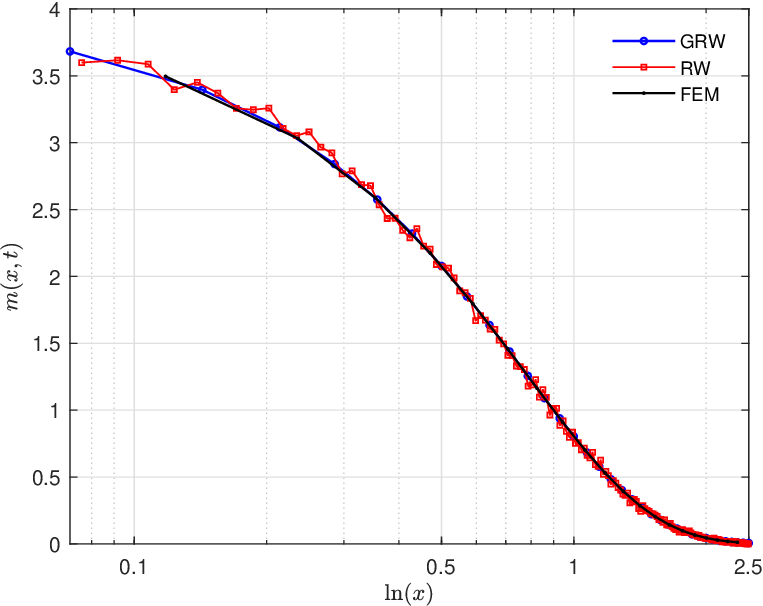}
%\caption{\label{fig:c_rubber}{\footnotesize \red{old} Diffusant profile as functions of space at time $T = 31$ minutes.}}
%\end{minipage}
\hspace{0.1cm}
%\begin{minipage}[c]{0.48\linewidth} %\centering
%\includegraphics[width=0.48\linewidth]{eps/s_rubber.eps}
\includegraphics[width=0.48\linewidth]{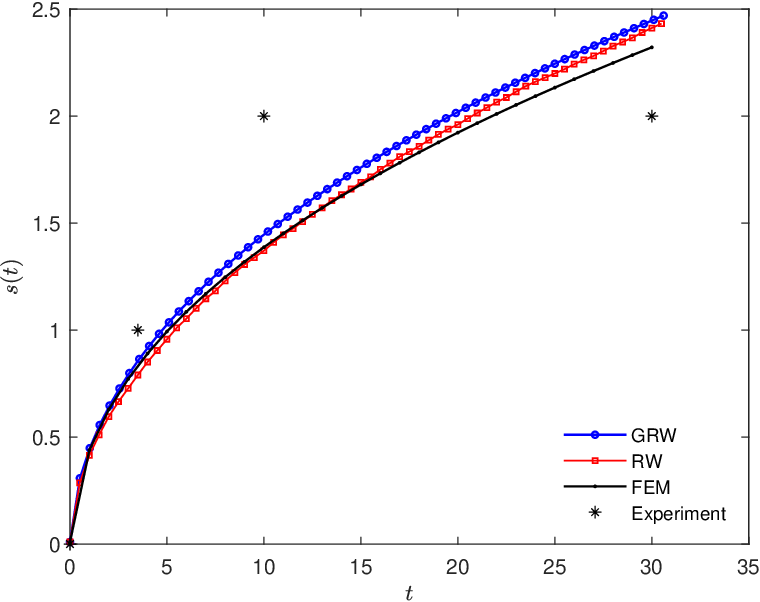}
\caption{\label{fig:s_rubber}{{\it Left}: Semilog plot for the profile of $m(x, t)$ as functions of space variable $x$ at time $T = 31$ minutes. {\it Right}: Comparison of the position of the moving front approximated by GRW, RW, and FEM with laboratory experimental data.}}
%\end{minipage}
\end{figure}

 The left-hand side of  Figure~\ref{fig:s_rubber} 
 shows the GRW simulation results compared to the RW and FEM solutions for the diffusant concentration profile. Similarly, the right-hand side of Figure~\ref{fig:s_rubber} compares the GRW solution with the RW and FEM solutions as well as with experimental data for the evolution of the penetration front. The GRW simulation results, for the concentration profile and the penetration fronts, are in good agreement with the RW and FEM results and are consistent with the laboratory experimental measurements of the penetration fronts.

 %\begin{remark}
While the RW concentration solution presented in Figure \ref{fig:s_rubber}, obtained with $500$ random walkers in \cite{Nepaletal2023}, is affected by oscillations, the GRW solution using $10^{24}$ particles is smooth and practically coincides with the FEM solution. We also note that the second-order convergence of the GRW scheme for diffusion in rubbers (see Table~\ref{tab:conv_rubber}) can only be obtained if the number of particles is at least $\mathcal{N}\approx 6.7\cdot 10^5$ (intuitively, for recovering much of the high-regularity of the solution).

Looking at Figure \ref{fig:s_rubber}, we notice that the numerically simulated penetration depths are within the expected physical range, i.e. they are closed to the measurement points. We refer the reader to  \cite{nepal2021moving} for the technical details of the experimental setup. Of course, one could perform a parameter identification procedure to identify the likelihood that our simulated moving boundary positions meet precisely the four measurements, but this is outside the scope of this work.
 %\end{remark}

%{\bf A TABLE with the order of convergence between GRW iterates (like a Cauchy sequence argument) is here missing... Please add}

\section{Conclusion}
\label{discussion}

In this work, we presented the GRW approach approximating numerically one-dimensional moving boundary problems. To the best of our knowledge, this is the first time when one-dimensional moving-boundary problems are handled by means of the GRW methodology.  

By means of numerical simulations, we propose a  GRW approximation and perform a convergence study for the classical Stefan problem and for a Stefan-type problem with a kinetic condition, with direct application to a more practical scenario related to the penetration of diffusants into rubber-based materials.   %We begin by exploring the GRW simulation results and comparing them with FEM simulation results for the classic Stefan problem, using its analytical solution to compute errors. 
For all treated cases,  both GRW and FEM approximations exhibited at least first-order of convergence for the position of the moving boundary and the concentration profile. 
%%(in suitable $L^2$-norms)
%\red{This is not accurate, we have different orders for the three problems investigated (see commented text below in the main.tex file).}
This is consistent with the existing convergence proofs for the FEM approximation. However, in practice, improved convergence rates may be attained, as indicated in the present study by the numerically inferred EOC for solutions of the Stefan-type problems with kinetic conditions.
This is often the case in numerical applications, see for instance~\cite{Radu2009}. 

Unlike the classical random walk method, the GRW method required less computing time than FEM. This makes one expect that, at least from the point of view of computing time,  the GRW method is suitable for moving-boundary problems posed in higher-dimensional settings. 

It remains an interesting and challenging task to prove rigorously uniform estimates on the indicated order of convergence. It is worth noting that remotely related work has been done very recently in \cite{Cuchiero} for the classical Stefan problem, the kinetic case remaining fully open.
%In the Stefan-type problem with a kinetic condition,  we achieved the first-order convergence for the position of the moving boundary, similar to the classic Stefan problem. However, the order of convergence for the solution to the corresponding partial differential equation increased to $1.5$. For the case of the moving boundary describing the penetration of a diffusant into rubber, we obtained a second order of convergence both for the position of the moving boundary and for the diffusant concentration profile.   

%We applied the GRW method to solve the Stefan problem with a kinetic condition, for which we constructed the analytical solution. In this scenario,  we achieve first-order convergence for the position of the moving boundary and $1.5$ order convergence for the solution to the corresponding partial differential equation. We next investigate the GRW method for a moving boundary problem describing the liquid penetration into rubbers.  In this case,  we constructed the manufacturing solution by adding source terms to the differential equation and the boundary conditions. This allowed us to compare the GRW solution to the exact solution and investigate the order of convergence in this case. We obtain the optimal orders of convergence, with a second order of convergence for both the position of the moving boundary and the diffusant concentration profile.   

\section*{Acknowledgments}
We thank U. Giese (Deutsches Institut f\"ur Kautschuktechnologie (DIK), Hannover, Germany) and N.H. Kr\"oger (tesa SE, formerly with DIK) for posing us this problem. Fruitful discussions with K. Kumazaki (Kyoto, Japan) and T. Aiki (Tokyo, Japan) on the mathematical analysis  of this moving-boundary problem are much appreciated.  The work of S.N.  and A.M. is partially funded by the Swedish Research Council's project {\em "Homogenization and dimension reduction of thin heterogeneous layers"} with grant nr. VR 2018-03648. A.M. acknowledges also support from the Knowledge Foundation for the grants KK 2019-0213, KK 2020-0152, and KK 2023-0010.

%----------------------------------------

\begin{thebibliography}{9}

\bibitem{Gupta2003}
{\sc S. C. Gupta}, {\it The Classical Stefan Problem: Basic Concepts, Modelling and Analysis}, Elsevier 2003.
\href{https://doi.org/10.1016/C2017-0-02306-6}{https://doi.org/10.1016/C2017-0-02306-6}

\bibitem{Stefan1891}
{\sc J. Stefan}, {\it \"{U}ber die Theorie der Eisbildung, insbesondere \"{u}ber
die Eisbildung im Polarmeere}, Ann. Physik Chemie {\bf 42} (1891), 269--286.
\href{https://doi.org/10.1002/andp.18912780206}{https://doi.org/10.1002/andp.18912780206}

\bibitem{tarzia1997one}
{\sc D. A Tarzia } and {\sc C. V. Turner},
{\it The one-phase supercooled Stefan problem with a convective boundary condition}, Q. Appl. Math. {\bf 55}(1) (1997), 41--50.
\href{ https://doi.org/10.1090/qam/1433750}{ https://doi.org/10.1090/qam/1433750}

\bibitem{Visintin}
{\sc A. Visintin}, {\it Models of {P}hase {T}ransition}, Birkh\"{a}user, 1996.

\bibitem{Tsunoda}
{\sc K. Tsunoda}, {\it Derivation of Stefan problem from a one-dimensional exclusion process with speed change}, Markov Process. Relat. {\bf 21} (2015), 263--273.

\bibitem{nepal2023analysis}
{\sc S. Nepal, Y. Wondmagegne} and {\sc  A. Muntean},
  {\it Analysis of a fully discrete approximation to a moving-boundary problem describing rubber exposed to diffusants}, Appl.  Math. Comput. {\bf 442} (2023), 127733.
  \href{https://doi.org/10.1016/j.amc.2022.127733}{https://doi.org/10.1016/j.amc.2022.127733}

\bibitem{nepal2021moving}
{\sc S. Nepal, R. Meyer, N. H. Kr{\"o}ger, T. Aiki, A. Muntean, Y. Wondmagegne} and {\sc U. Giese},
  { \it A moving boundary approach of capturing diffusants penetration into rubber: {FEM} Approximation and comparison with laboratory measurements}, KGK-Kaut. Gumi. Kunst. {\bf 5} (2021), 61--69.
  
\bibitem{Villa}
{\sc D. A. Tarzia,  L. T. Villa}, {\it On the free boundary problem in the Wen-Langmuir shrinking core model for noncatalytic gas-solid reactions}, Meccanica {\bf 24} (1989), 86--92. \href{https://doi.org/10.1007/BF01560134}{https://doi.org/10.1007/BF01560134}

\bibitem{Aiki}
{\sc T. Aiki,  A. Muntean}, {\it Existence and uniqueness of solutions to a matheatical model predicting service life of concrete structures}, Advances in Mathematical Sciences and Applications {\bf 19}(1) (2009), 109--129. 

\bibitem{Evans}
{\sc J. D. Evans,   J. R. King},
{\it The Stefan problem with nonlinear kinetic undercooling},
Q. J. Mech. Appl. {\bf 56}(1) (2003), 139--161. \href{https://doi.org/10.1093/qjmam/56.1.139}{https://doi.org/10.1093/qjmam/56.1.139}

\bibitem{Nepaletal2023}
{\sc S. Nepal, M. \"{O}gren, Y. Wondmagegne} and {\sc A. Muntean}, {\it Random walks and moving boundaries: Estimating the penetration of diffusants into dense rubbers}, Probabilist. Eng. Mech. {\bf 74} (2023), 103546.
\href{https://doi.org/10.1016/j.probengmech.2023.103546}{https://doi.org/10.1016/j.probengmech.2023.103546}

\bibitem{Suciuetal2021}
{\sc N. Suciu, D. Illiano, A. Prechtel} and {\sc F.A. Radu}, {\it Global random walk solvers for fully coupled flow and transport in saturated/unsaturated porous media}, Adv. Water Resour. {\bf 152} (2021), 103935.
\href{https://doi.org/10.1016/j.advwatres.2021.103935}{https://doi.org/10.1016/j.advwatres.2021.103935}

\bibitem{Suciu2019}
{\sc N. Suciu}, {\it Diffusion in Random Fields. Applications to Transport in Groundwater}, Birkh\"{a}user, Cham, 2019.
\href{https://doi.org/10.1007/978-3-030-15081-5}{https://doi.org/10.1007/978-3-030-15081-5}

\bibitem{SuciuandRadu2022}
{\sc N. Suciu} and {\sc F.A. Radu}, {\it Global random walk solvers for reactive transport and biodegradation processes in heterogeneous porous media}, Adv. Water Resour. {\bf 166} (2022), 104268.
\href{https://doi.org/10.1016/j.advwatres.2022.104268}{https://doi.org/10.1016/j.advwatres.2022.104268}

\bibitem{Strikwerda2004}
{\sc J.C. Strikwerda}, {\it Finite Difference Schemes and Partial Differential Equations}, SIAM, 2004.
\href{https://doi.org/10.1137/1.9780898717938}{https://doi.org/10.1137/1.9780898717938}

\bibitem{Suciuetal2024}
{\sc N. Suciu, F.A. Radu} and {\sc E. C\u{a}tina\c{s}}, {\it Iterative schemes for coupled flow and transport in porous media -- Convergence and truncation errors}, Numer. Anal. Approx. Theory {\bf 53}(1) (2024), 158--183.
\href{https://doi.org/10.33993/jnaat531-1429}{https://doi.org/10.33993/jnaat531-1429}

\bibitem{kutluay1997numerical}
{\sc S. Kutluay, A.R. Bahadir} and {\sc A. {\"O}zde{\c{s}}},
{\it The numerical solution of one-phase classical Stefan problem}, J. Comput. Appl. Math. {\bf 81}(1) (1997), 135--144. \href{https://doi.org/10.1016/S0377-0427(97)00034-4}{https://doi.org/10.1016/S0377-0427(97)00034-4}

\bibitem{savovic2003finite}
{\sc S. Savovi{\'c}} and {\sc J. Caldwell},
{\it Finite difference solution of one-dimensional Stefan problem with periodic boundary conditions}, Int. J. Heat. Mass Tran. {\bf 46}(15) (2003), 2911--2916. \href{https://doi.org/10.1016/S0017-9310(03)00050-4}{https://doi.org/10.1016/S0017-9310(03)00050-4}

\bibitem{mori1977finite}
{\sc M. Mori},
{\it A finite element method for solving the two phase Stefan problem in one space dimension}, Publ. Res. I. Math. Sci. {\bf 13}(3) (1977), 723--753. \href{https://doi.org/10.2977/prims/1195189605}{https://doi.org/10.2977/prims/1195189605}

\bibitem{mori1976stability}
{\sc M. Mori},
{\it Stability and convergence of a finite element method for solving the Stefan problem}, Publ. Res. I. Math. Sci. {\bf 12}(2) (1976), 539--563. \href{https://doi.org/10.2977/prims/1195190728}{https://doi.org/10.2977/prims/1195190728}



\bibitem{casaban2023numerical}
  {\sc M.-C. Casab{\'a}n, R. Company} and {\sc L. J{\'o}dar},
  {\it Numerical difference solution of moving boundary random Stefan problems},
  Math. Comput. Simulat.
  {\bf 205} (2023), 878--901. \href{https://doi.org/10.1016/j.matcom.2022.10.026}{https://doi.org/10.1016/j.matcom.2022.10.026}


  \bibitem{Ogren2022}
{\sc M. \"{O}gren}, {\it Stochastic solutions of Stefan problems with general time-dependent
boundary conditions}, in: A. Malyarenko, Y. Ni, M. Ran\v{c}i\'{e}, S. Silvestrov
(Eds.), Stochastic Processes, Statistical Methods, and Engineering Mathematics.
SPAS 2019, Springer Proceedings in Mathematics \& Statistics, vol.
{\bf 408}, Springer, Cham, 2022.
\href{http://dx.doi.org/10.1007/978-3-031-17820-7_29}{http://dx.doi.org/10.1007/978-3-031-17820-7\_29}
% (arXiv:2006.04939)

\bibitem{Roy2005}
{\sc C.J. Roy}, {\it Review of code and solution verification procedures for computational simulation}, J. Comput. Phys. {\bf 205} (2005), 131--156. \href{https://doi.org/10.1016/j.jcp.2004.10.036}{https://doi.org/10.1016/j.jcp.2004.10.036}

\bibitem{Alecsaetal2020}
{\sc C.D.Alecsa, I. Boros, F. Frank, P. Knabner, M. Nechita, A. Prechtel, A. Rupp} and {\sc N. Suciu}, {\it Numerical benchmark study for flow in heterogeneous aquifers}, Adv. Water Resour. {\bf 138} (2020), 103558.
\href{https://doi.org/10.1016/j.advwatres.2020.103558}{https://doi.org/10.1016/j.advwatres.2020.103558}

\bibitem{lewis1980laboratory}
{\sc P. M. Lewis},
  {\it Laboratory testing of rubber durability},
  Polom. Test. {\bf 1} (1980), 167--189. \href{https://doi.org/10.1016/0142-9418(80)90002-1}{https://doi.org/10.1016/0142-9418(80)90002-1}

 \bibitem{hayes1955diffusion} 
  {\sc M.J. Hayes} and {\sc G.S. Park},
  {\it The diffusion of benzene in rubber. Part 1.—Low concentrations of benzene},
  T. Faraday Soc. {\bf 51} (1955), 1134--1142. \href{https://doi.org/10.1039/TF9555101134}{https://doi.org/10.1039/TF9555101134}

  
\bibitem{kaliyathan2021analysis}
{\sc A. V. Kaliyathan, A. V. Rane, S. Jackson} and {\sc S. Thomas},
  {\it Analysis of diffusion characteristics for aromatic solvents through carbon black filled natural rubber/butadiene rubber blends}, Polym. Composite {\bf 42} (2021), 375--396. \href{https://doi.org/10.1002/pc.25832}{https://doi.org/10.1002/pc.25832}

 \bibitem{kumazaki2020global}
   {\sc  K. Kumazaki} and {\sc A. Muntean},
  {\it Global weak solvability, continuous dependence on data, and large time growth of swelling moving interfaces}, Interf. Free Bound. {\bf 22}(1) (2020), 27--50. \href{https://doi.org/10.4171/ifb/431}{https://doi.org/10.4171/ifb/431}

 \bibitem{nepal2021error}
{\sc S. Nepal, Y. Wondmagegne} and {\sc  A. Muntean},
  {\it Error estimates for semi-discrete finite element approximations for a moving boundary problem capturing the penetration of diffusants into rubber}, Int. J. Numer. Anal. Mod. {\bf 19} (2022), 101--125. \href{https://global-sci.org/intro/article_detail/ijnam/20351.html}{https://global-sci.org/intro/article\_detail/ijnam/20351.html}

  \bibitem{Roache2002}
{\sc P. J. Roache}, {\it Code verification by the method of manufactured solutions}, J. Fluids Eng. {\bf 124}(1) (2002), 4--10.
\href{http://dx.doi.org/10.1115/1.1436090}{http://dx.doi.org/10.1115/1.1436090}

\bibitem{Cuchiero}
{\sc C. Cuchiero, C. Reisinger} and {\sc S. Rigger}, {\it Implicit and fully discrete approximation of the supercooled Stefan problem in the presence of blow-ups}, SIAM J. Numer. Anal. {\bf 62}(3) (2024), 1145--1170. \href{https://doi.org/10.1137/22M1509722}{https://doi.org/10.1137/22M1509722}

\bibitem{Radu2009} 
{\sc F. A. Radu, I. S. Pop} and {\sc S. Attinger}, {\it Analysis of an Euler implicit-mixed finite element scheme for reactive solute transport in porous media}, Numer. Meth. Part. D. E. {\bf 26}(2) (2009), 320--344.
\href{https://doi.org/10.1002/num.20436}{https://doi.org/10.1002/num.20436} 


\end{thebibliography}
\end{document}